# ON POSITIVE INTEGERS REPRESENTED AS ARITHMETIC SERIES

**Dimitris Sardelis**

**OVERVIEW** The aim of the present article is to explore the possibilities of representing positive integers as sums of other positive integers and highlight certain fundamental connections between their multiplicative and additive properties. In particular, we shall be concerned with the representation of positive integers as arithmetic series of the simplest kind, i.e., either as sums of successive odd positive numbers, or as sums of successive even positive numbers (both treated as Problem 1), or as sums of consecutive positive integers (treated as Problem 2).

## INTRODUCTION

It is a historic fact that generating natural numbers from arithmetical progressions and representing them with geometrical forms, originated and shaped invariably early number theory (e.g. see [**1**] and [**2**]). The formation of triangular numbers $1, 1+2, 1+2+3, \ldots$, and of square numbers $1, 1+3, 1+3+5, \ldots$, are the most well known examples. Less known perhaps is that Diophantus [**3**] defines polygonal numbers as arithmetic series and even less known that Nicomachus [**4**] similarly defines three dimensional pyramidal numbers. Extending this trend further, one may produce number formations of higher dimensionality [**5**].

The present exposition is intended to follow the opposite direction: Dicossiating number and geometrical form, our aim is to tackle the inverse problem of the possible ways that any positive integer can be expressed as an arithmetic series of other positive integers. In what follows we restrict ourselves to simple arithmetic series representations, i.e., either sums of successive odd/even positive numbers (Problem 1) or sums of consecutive positive numbers (Problem 2) . The prime factorization form of numbers proves to play a key role in such representations.

**PROBLEM 1** ***It is required to represent every positive integer ≻ 1 either by sums of successive odd positive numbers or by sums of successive even positive numbers.***

Let $a, a+2, a+4, \ldots, b$ be an arithmetic progression of $r$ positive integers with a common difference of $2$. Then any other positive integer $N$ can be expressed by the series $a+(a+2)+(a+4)+\cdots+b$ if $N = r\cdot(r+a-1).$ On the other hand, it is clear that every integer $>1$ can always be written as the product of any two of its complementary divisors, i.e., as $N = d\cdot d'$ where d is any of the divisors of $N.$ Hence, we must have: $d\cdot d' = r\cdot(r+a-1).$ Since $a\geq 1$, then $r+a-1\geq r.$ Consequently, $d = r$ and $d' = r+a-1$ for each $d\leq d'.$ Therefore, for every pair of complementary divisors $d, d'$ of $N$ such that $d\leq d'$, the first and last terms of the series are $a = d'-d+1$ and $b = d'+d-1$, respectively.

- If $N$ is prime, then it has a single pair of complementary divisors: $d = 1, d' = N.$ Therefore, $a = b = N$ and $N$ is representable neither by sums of successive odd numbers nor by sums of successive even numbers.
- If *N* is composite, then it has at least one divisor $d>1$ and, subsequently, it can be represented by at least one of the aforementioned sums. These conclusions may be formally stated as follows:

**THEOREM** 1′ *No prime is expressible as a sum of successive odd positive numbers or as a sum of successive even positive numbers.*

**THEOREM 1** *Every positive integer $N$ with at least one divisor $d>1,$ is expressible as an arithmetic series of $d$ positive integers with a common difference of $2$. Thus, for each pair of complementary divisors $d$ and $d'$ of $N$ such that $1<d\leq d',$ $N$ can be written as*

$$N = (d'-d+1)+(d'-d+3)+\cdots+(d'+d-1). \qquad (1)$$

It is clear that if $d'-d$ is even, sum (1) consists of $d$ successive odd positive numbers while if $d'-d$ is odd, it consists of $d$ successive even positive numbers. The following two Corollaries may thus be stated:

**Corollary 1.1** *If any two complementary divisors $d, d'$ of a positive integer $N$ differ by an even number, i.e., if $d'-d = 2\cdot(c-1)$ with $c$ a positive integer, then $N$ can be expressed as a sum of $d$ successive odd positive numbers:*

$$\begin{gathered} N = d\cdot[d+2\cdot(c-1)] = \\ [2\cdot(c-1)+1]+[2\cdot(c-1)+3]+\cdots+[2\cdot(c-1)+2d-1]. \end{gathered} \qquad (2)$$

As examples, the number $120$ can be written as $10\cdot 12 = 3+5+\cdots+21$ ($10$ terms), and the number $187$ can be written as $11\cdot 17 = 7+9+11+\cdots+27$ ($11$ terms).

In particular, for $d = d'$ $[c = 1]$ we have the Proposition:

**Proposition 1.1** *Every squared number $N = d^2$ can be written as the sum of the first $d$ odd positive numbers:*

$$N = d^2 = 1+3+5+\cdots+(2d-1). \qquad (3)$$

For $d = 2, 3, 4, \ldots, 9,$ we have Table 1:

| Number | Sum | # of terms |
|---|---|---|
| $2^2$ | 1 + 3 | 2 |
| $3^2$ | 1 + 3 + 5 | 3 |
| $4^2$ | 1 + 3 + 5 + 7 | 4 |
| $5^2$ | 1 + 3 + 5 + 7 + 9 | 5 |
| $6^2$ | 1 + 3 + 5 + 7 + 9 + 11 | 6 |
| $7^2$ | 1 + 3 + 5 + 7 + 9 + 11 + 13 | 7 |
| $8^2$ | 1 + 3 + 5 + 7 + 9 + 11 + 13 + 15 | 8 |
| $9^2$ | 1 + 3 + 5 + 7 + 9 + 11 + 13 + 15 + 17 | 9 |

**Corollary 1.2** *If any two complementary divisors* $d, d'$ *of a positive integer* $N$ *differ by an odd number, i.e., if* $d' - d = 2c - 1$ *with* $c$ *a positive integer, then* $N$ *can be expressed as a sum of* $d$ *successive even positive numbers:*

$$N = d\cdot(d + 2c - 1) \\ = 2c + 2\cdot(c+1) + 2\cdot(c+2) + \cdots + 2\cdot(c+d-1). \tag{4}$$

As examples, the number $150$ can be written as $10\cdot 15 = 6 + 8 + 10 + \cdots + 24$ ($10$ terms), and the number $198$ can be written as $11\cdot 18 = 8 + 10 + 12 + \cdots + 28$ ($11$ terms).

In particular, for $d' = d + 1\,[c = 1]$, we have the Proposition:

**Proposition 1.2** *Every positive integer of the form* $N = d\cdot(d+1)$ *can be written as the sum of the first* $d$ *even positive numbers:*

$$N = d\cdot(d+1) = 2 + 4 + 6 + \cdots + 2d. \tag{5}$$

For $d = 2, 3, 4, \ldots, 9,$ we have Table 2:

| Number | Sum | # of terms |
|---|---|---|
| 2 · 3 | 2 + 4 | 2 |
| 3 · 4 | 2 + 4 + 6 | 3 |
| 4 · 5 | 2 + 4 + 6 + 8 | 4 |
| 5 · 6 | 2 + 4 + 6 + 8 + 10 | 5 |
| 6 · 7 | 2 + 4 + 6 + 8 + 10 + 12 | 6 |
| 7 · 8 | 2 + 4 + 6 + 8 + 10 + 12 + 14 | 7 |
| 8 · 9 | 2 + 4 + 6 + 8 + 10 + 12 + 14 + 16 | 8 |
| 9 · 10 | 2 + 4 + 6 + 8 + 10 + 12 + 14 + 16 + 18 | 9 |

**Table 2**

Theorem 1 provides the framework for expressing any given positive integer by sums of successive odd or successive even positive numbers. In fact, all such series representations of a given positive integer can be determined by using Theorem 1. On the other hand, Corollaries 1.1 and 1.2 specify those structural elements in the factorization of numbers which determine whether the summands in their series representations are all odd or all even. A question then arises: Are there positive integers which are exclusively representable either by sums of successive odd or by sums of successive even numbers? The answer will be articulated by stating a few propositions that apply to particular number forms.

It is most natural to start with the series representations of powers of single primes where there is a clear-cut distinction and ordering of the pairs of complementary divisors. Thus, we have the following Proposition:

**Proposition 1.3** *Numbers of the form $p^{n+1}$, where $p$ is any prime and $n$ is a positive integer, can be exclusively expressed by $\lfloor (n+1)/2 \rfloor$ sums of $p^k$ successive odd positive numbers, where $k = 1, 2, 3, \ldots, \lfloor (n+1)/2 \rfloor$(The symbol $\lfloor x \rfloor$ stands for the largest integer not exceeding $x$ ). Thus, for each $k$, we have*

$$p^{n+1} = \left(p^{n-k+1} - p^k + 1\right) + \left(p^{n-k+1} - p^k + 3\right) + \cdots + \left(p^{n-k+1} + p^k - 1\right). \qquad (6)$$

As examples, the number $2^4$ is expressed by the sums $7 + 9$ ($2$ terms), $1 + 3 + 5 + 7$ ($2^2$ terms), and the number $3^5$ is expressed by the sums $79 + 81 + 83$ ($3$ terms), $19 + 21 + 23 + \cdots + 35$ ($3^2$ terms).

Particular cases of Proposition 1.3 are the (sub)propositions:

**Proposition 1.3a** *No prime raised to an odd power is expressible either as a sum of successive odd positive numbers starting from $1$, or as a sum of successive even positive numbers starting from $2$.*

Every summand in (6) is odd. Therefore, each series representation of $p^{n+1}$ is an odd series. Furthermore, if $n+1$ is odd, then $2k < n+1$ and $p^{n-k+1} > p^k$. Hence, the first term in each series is $> 2$.

**Proposition 1.3b** Every *prime raised to an even power can be written as a sum of successive odd positive numbers starting from $1$.*

This follows readily from Proposition 1.2 by setting $n \to 2n$ and $k = n$.

The next step is to consider the series representations of composite numbers. For the case of composite odd numbers we have the following proposition:

**Proposition 1.4** *Every composite odd positive number is expressible exclusively by sums of successive odd positive numbers.*

Let a composite odd positive number have the prime factorization form $p_1 p_2 \cdots p_L$, where $p_1, p_2, p_3, \ldots,$ are odd primes, not necessarily distinct: $p_1 \leq p_2 \leq p_3 \leq \cdots \leq p_L$. Since all the divisors of this number are odd, any two of its complementary divisors $d, d'$ differ by an even number and all summands in its series representations (1) are odd. Furthermore, among all possible pairs $d, d'$ there is at least one for which $1 < d \leq d'$. Consequently, the composite odd numbers are representable by at least one sum of successive odd positive numbers.

We illustrate Proposition 1.4 by listing in Table 3 all series representations of the composite numbers $p_1 p_2$, $(p_1 p_2)^2$ and $(p_1 p_2)^3$ as sums of successive odd positive integers, where $p_1, p_2$ are any odd primes so that $p_1 < p_2$.

| Number | Series Representation | # of terms |
|---|---|---|
| $p_1 p_2$ | $(p_2 - p_1 + 1) + (p_2 - p_1 + 3) + \cdots + (p_2 + p_1 - 1)$ | $p_1$ |
| $(p_1 p_2)^2$ | $(p_1 p_2^2 - p_1 + 1) + (p_1 p_2^2 - p_1 + 3) + \cdots + (p_1 p_2^2 + p_1 - 1)$ | $p_1^1$ |
| | $(p_1^2 p_2 - p_2 + 1) + (p_1^2 p_2 - p_2 + 3) + \cdots + (p_1^2 p_2 + p_2 - 1)$ | $p_2^1$ |
| | $(p_2^2 - p_1^2 + 1) + (p_2^2 - p_1^2 + 3) + \cdots + (p_2^2 + p_1^2 - 1)$ | $p_1^2$ |
| | $1 + 3 + 5 + 7 + \cdots + (2 p_1 p_2 - 1)$ | $p_1^1 p_2^1$ |
| $(p_1 p_2)^3$ | $(p_1^2 p_2^3 - p_1 + 1) + (p_1^2 p_2^3 - p_1 + 3) + \cdots + (p_1^2 p_2^3 + p_1 - 1)$ | $p_1^1$ |
| | $(p_1^3 p_2^2 - p_2 + 1) + (p_1^3 p_2^2 - p_2 + 3) + \cdots + (p_1^3 p_2^2 + p_2 - 1)$ | $p_2^1$ |
| | $(p_1 p_2^3 - p_1^2 + 1) + (p_1 p_2^3 - p_1^2 + 3) + \cdots + (p_1 p_2^3 + p_1^2 - 1)$ | $p_1^2$ |
| | $(p_1^2 p_2^2 - p_1 p_2 + 1) + (p_1^2 p_2^2 - p_1 p_2 + 3) + \cdots + (p_1^2 p_2^2 + p_1 p_2 - 1)$ | $p_1^1 p_2^1$ |
| | $(p_1 p_2^2 - p_1^2 p_2 + 1) + (p_1 p_2^2 - p_1^2 p_2 + 3) + \cdots + (p_1 p_2^2 + p_1^2 p_2 - 1)$ | $p_1^2 p_2$ |
| | $(p_1^3 - p_2^3 + 1) + (p_1^3 - p_2^3 + 1) + \cdots + (p_1^3 + p_2^3 - 1)$ | $p_1^3$ |
| | $(p_2^2 - p_1^3 p_2 + 1) + (p_2^2 - p_1^3 p_2 + 1) + \cdots + (p_2^2 + p_1^3 p_2 - 1)$ | $p_1^3 p_2$ if $p_1^3 < p_2$ |
| | $(p_1^3 p_2 - p_2^2 + 1) + (p_1^3 p_2 - p_2^2 + 3) + \cdots + (p_1^3 p_2 + p_2^2 - 1)$ | $p_2^2$ if $p_1^2 > p_2$ |

**Table 3**

It is interesting to note that the series representation of the number $(p_1 p_2)^2$ which consists of $p_1 p_2$ terms, can also be written as

$$N^2 = 1 + 3 + 5 + \cdots + (2 N - 1),$$

with $N = p_1 p_2.$ Similarly, the series representation of number $(p_1 p_2)^3$ which consists of $p_1 p_2$ terms, can also be written as

$$N^3 = \left(N^2 - N + 1\right) + \left(N^2 - N + 3\right) + \cdots + \left(N^2 + N - 1\right).$$

Thus, the series representations of composite odd numbers raised to powers $> 1$, do not depend on their prime factorization form.This observation is suggestive of the following Proposition:

**Proposition 1.**4′ *Numbers of the form $N^{n+1}$, where $N$ is odd $> 1$ and $n$ is a positive integer, can be expressed – not exclusively – by $\lfloor (n + 1) / 2 \rfloor$ sums of $N^k$ successive odd positive numbers.Thus, for each $k = 1, 2, 3, \ldots, \lfloor (n + 1) / 2 \rfloor$, we have*

$$N^{n+1} = \left(N^{n-k+1} - N^k + 1\right) + \left(N^{n-k+1} - N^k + 3\right) + \cdots + \left(N^{n-k+1} + N^k - 1\right). \qquad (7)$$

Numbers of the form $N^{n+1}$ where $N$ is odd, can always be factored so that $d = N^k$, $d' = N^{n-k+1}$ with $2 k \leq n + 1.$ Hence, Theorem 1 readily applies to yield (7). It is obvious that if $N$ is prime, (7) reduces to (6).

According to whether the exponent of $N$ is odd or even, we have the following particular Propositions:

**Proposition 1.**4′ a *No number of the form $N^{2n+1}$, where $N$ is odd $> 1$ and $n$ is a positive integer, can be written as a sum of successive odd positive numbers starting from $1.$*

Setting $n \to 2 n$ in (7), it is clear that $2 k < 2 n + 1.$ Consequently, each summand in the series representation of $N^{2n+1}$ exceeds $1.$

**Proposition 1.**4′ b *Every number of the form $N^{2n}$, where $N$ is odd $> 1$ and n is a positive integer, can be written as*

*the sum of the first $N^n$ odd positive numbers.*

Thus, one recovers Proposition 1.1 in a restated form. Proposition 1.4′ b is deduced by setting $n+1 \to 2n$ and $k = n$ in (7).

As a numerical illustration of Propositions 1.4′ a and 1.4′ b, the number $15^4$ has twelve divisors $d$ such that $1 < d \leq d'$:

$$3,\ 5,\ 3^2,\ 3\cdot 5,\ 5^2,\ 3^3,\ 3^2\cdot 5,\ 3\cdot 5^2,\ 3^4,\ 5^3,\ 3^3\cdot 5,\ 3^2\cdot 5^2,$$

and is expressible by twelve sums of successive odd positive numbers. The prime factorization of $15$ does not enter in two of these sums ($d = 3\cdot 5,\ 3^2\cdot 5^2$). For the case $d = 3^2\cdot 5^2$, the number $15^4$ can be written as the sum of the first $15^2$ odd positive numbers [Proposition 1.4′ a]. On the other hand, the number $15^5$ can also be expressed by two sums not depending on its prime factorization (for $d = 3\cdot 5,\ 3^2\cdot 5^2$), but neither of them starts from $1$ [Proposition 1.4′ b].

We have thus reached a partial answer to the question posed above. Our findings can be summarized as follows: Composite odd positive numbers and even numbers of the form $2^{n+1}$, where $n$ is a positive integer, are exclusively representable by sums of successive odd positive numbers. Stated somewhat differently, products of even primes (powers of $2$ greater than $1$) and products of odd primes, are each exclusively expressible by sums of successive odd positive numbers.

The second version of this summary provides the insight for a complete answer: According to the fundamental theorem of arithmetic, any number can be written uniquely as a product of primes, not necessarily distinct. Consequently, expressible exclusively by sums of successive odd positive numbers, are all the positive integers for which every two complementary divisors differ by an even number [Corollary 1.1]. Therefore, the divisors of these integers must be either all even or all odd. Hence,the following Theorem can be stated:

**THEOREM 1A** *Expressed exclusively by sums of successive odd positive numbers are all the composite odd positive numbers and all the even numbers of the form $2^{n+1}$, with n a positive integer .*

On the other hand, expressible exclusively by sums of successive even positive numbers are all the positive integers for which the difference of any two of its complementary divisors is odd [ Corollary 1.2 ]. Therefore, the divisors of these integers cannot be all odd and no more than one of them can be even. Consequently, the prime factorization of the integers expressible exclusively by sums of successive even positive numbers, must consist of a single even prime and a product of odd primes. Hence,the following Theorem can be stated:

**THEOREM 1B** *Expressed exclusively by sums of successive even positive numbers are all the even numbers of the form $2\cdot(2m+1)$, where $m$ is a positive integer.*

As a numerical example of Theorem 1.B, the number $90$ can be expressed exclusively by the following sums of successive even positive numbers: $44+46$ ($2$ terms), $28+30+32$ ($3$ terms), $14+16+18+20+22$ ($5$ terms), $10+12+14+\cdots+20$ ($2\cdot 3$ terms), and $2+4+6+\cdots+18$ ($3^2$ terms).

By Proposition 1.2, a number of the form $2\cdot(2m+1)$ with a pair of complementary divisors $d$ and $d'$ so that $d' = d+1$, can be written as a sum of the first $d$ even positive numbers.The next two Propositions refer to such representations.

**Proposition 1.5a** *Every integer of the form $2\cdot(2m+1)$ for which $m = k\cdot(4k-3)$, where $k$ is a positive*

*integer, can be written as the sum of the first $4k-2$ even positive numbers.*

Thus, for $k = 1, 2, 3, 4, 5,$ we have Table 4:

| k | Number | Sum | # of terms |
|---|---|---|---|
| 1 | $2\cdot 3$ | $2+4$ | 2 |
| 2 | $6\cdot 7$ | $2+4+\cdots+12$ | 6 |
| 3 | $10\cdot 11$ | $2+4+\cdots+20$ | 10 |
| 4 | $14\cdot 15$ | $2+4+\cdots+28$ | 14 |
| 5 | $18\cdot 19$ | $2+4+\cdots+36$ | 18 |

**Table 4**

**Proposition 1.5b** *Every integer of the form $2\cdot(2m+1)$ for which $m = k\cdot(4k+3)$, where $k$ is a positive integer, can be written as the sum of the first $4k+1$ even positive numbers.*

Thus, for $k = 1, 2, 3, 4, 5,$ we have Table 5:

| k | Number | Sum | # of terms |
|---|---|---|---|
| 1 | $5\cdot 6$ | $2+4+\cdots+10$ | 5 |
| 2 | $9\cdot 10$ | $2+4+\cdots+18$ | 9 |
| 3 | $13\cdot 14$ | $2+4+\cdots+26$ | 13 |
| 4 | $17\cdot 18$ | $2+4+\cdots+34$ | 17 |
| 5 | $21\cdot 22$ | $2+4+\cdots+42$ | 21 |

**Table 5**

Propositions 1.5a and 1.5b are deduced by requiring that $2\cdot(2m+1) = d\cdot(d+1)$. It then follows that the discriminant of this quadratic equation for $d$ must be the square of an odd number. Hence one finds that $1+8\cdot(2m+1) = (8k\pm 3)^2$. Consequently, $d = 4k-2$ or $d = 4k+1$ and, by Proposition 1.2, the corresponding sums are

$$2+4+6+\cdots+(8k-4) \quad and \quad 2+4+6+\cdots+(8k+2),$$

respectively.

Since $2$ divides all numbers of the form $2\cdot(2m+1)$, it is clear that

**Proposition 1.6** *Every number of the form $2\cdot(2m+1)$,where $m$ is a positive integer, is expressible as the sum of two successive even positive numbers.*

Obviously, the sum is $2m+2\cdot(m+1)$. If $m$ is a prime, then this is the only representation. Thus, we may state:

**Proposition 1.**6′ *Every number of the form $2p$, where $p$ is an odd prime, is expressible as the sum of no more than two successive even positive numbers.*

The only natural numbers for which Theorems 1A and 1B do not apply, are the even numbers of the form $2^{n+1}\cdot(2m+1)$, where $m$ and $n$ are positive integers. These numbers include divisors which, combined in complementary pairs, they generate both kinds of sums, i.e., sums of successive odd positive numbers as well as sums of successive even positive numbers. This third "mixed component" of Theorem 1, may be stated as follows:

**THEOREM 1C** *Expressed by sums of successive odd positive numbers as well as by sums of successive even positive numbers, are all even numbers of the form $2^{n+1}\cdot(2m+1)$, where $m$ and $n$ are positive integers.*

Indeed, a possible class of complementary divisors is $d = 2^k$ and $d' = 2^{n-k+1}$, where $k$ is a positive integer such that $2k < n+1$. Then for each $k$, the resulting series consists of successive odd positive numbers. A different class of complementary divisors is $d = 2^{n+1} < d' = 2m+1$, or $d = 2m+1 < d' = 2^{n+1}$. In either case, the resulting series consists of successive even positive numbers.

As a numerical example, the number $60$ can be expressed as a series of successive odd numbers by the two sums: $29+31$($2$ terms), $5+7+9+\cdots+15$ ($2 \cdot 3$ terms), and as a series of successive even numbers by the three sums: $18+20+22$ ($3$ terms), $12+14+16+18$ ($2^2$ terms), $8+10+12+14+16$ ($5$ terms).

By Propositions 1.5a and 1.5b, the sum of the first $4k-2$ or $4k+1$ even positive numbers yields even numbers of the form described by Theorem 1B. Since no power of $2$ can be expressed as a sum of successive even positive numbers starting from $2$ [Proposition 1.3a], it follows that the sum of the first $4k$ or $4k-1$ even positive numbers must yield even numbers of the form described by Theorem 1C. Thus we have the following Propositions:

**Proposition 1.7a** *Numbers of the form $2^2 k \cdot (4k-1)$, where $k$ is a positive integer, can be written as the sum of the first $4k-1$ even positive numbers.*

For $k = 1, 2, 3, 4, 5,$ we have Table 6:

| k | Number | Sum | # of terms |
|---|---|---|---|
| 1 | 3 · 4 | 2 + 4 + 6 | 3 |
| 2 | 7 · 8 | 2 + 4 + 6 + ⋯ + 14 | 7 |
| 3 | 11 · 12 | 2 + 4 + 6 + ⋯ + 22 | 11 |
| 4 | 15 · 16 | 2 + 4 + 6 + ⋯ + 30 | 15 |
| 5 | 19 · 20 | 2 + 4 + 6 + ⋯ + 38 | 19 |

**Table 6**

**Proposition 1.7b** *Numbers of the form $2^2 k \cdot (4k+1)$, where $k$ is a positive integer, can be written as the sum of the first $4k$ even positive numbers.*

For $k = 1, 2, 3, 4, 5,$ we have Table 7:

| k | Number | Sum | # of terms |
|---|---|---|---|
| 1 | 4 · 5 | 2 + 4 + 6 + 8 | 4 |
| 2 | 8 · 9 | 2 + 4 + 6 + ⋯ + 16 | 8 |
| 3 | 12 · 13 | 2 + 4 + 6 + ⋯ + 24 | 12 |
| 4 | 16 · 17 | 2 + 4 + 6 + ⋯ + 32 | 16 |
| 5 | 20 · 21 | 2 + 4 + 6 + ⋯ + 40 | 20 |

**Table 7**

Propositions 1.7a and 1.7b can be deduced by requiring that

$$2^{n+1} \cdot (2m+1) = d \cdot (d+1).$$

Since these numbers are divisible by $4$, either we must have that $d = 4k,\ d' = 4k-1$, or that $d = 4k-1,\ d' = 4k$, with $k$ a positive integer. Hence, Proposition 1.2 yields

$$4k \cdot (4k+1) = 2+4+6+\cdots+8k,$$

and

$$(4k-1)\cdot 4k = 2+4+6+\cdots+(8k-2).$$

For the pair of divisors $d=2$ and $d'=2^n\cdot(2m+1)$, we have that the integers $N=2^{n+1}\cdot(2m+1)$ can be written as the sum of two successive odd positive numbers: $[(N/2)-1]+[(N/2)+1]$. Stated in a proposition form:

**Proposition 1.8** *Every even number of the form $N=2^{n+1}\cdot(2m+1)$, where $m$ and $n$ are positive integers, can be written as the sum of two successive odd positive numbers.*

In particular, if $2m+1$ is an odd prime, then this is the only representation for N as a sum of successive odd numbers. Stated in a proposition form:

**Proposition** 1.8′ *Every even number of the form $2^2\,p$, where $p$ is an odd prime, can be expressed as the sum of no more than two successive odd positive numbers.*

The numbers of the form $2^{n+1}\cdot(2m+1)$ also have as complementary divisors the pair: $d=2^2,\ d'=2^{n-1}\cdot(2m+1)$. For the corresponding series representations we have the following three cases:

- For $m=n=1$, the corresponding number $2^2\cdot 3$ can be written as the sum of three successive even numbers: $2+4+6$.
- For $n=1$ and $m>1$, the numbers $2^2\cdot(2m+1)$ can be written as the sum of four successive even positive numbers:

$$2(m-1)+2m+2(m+1)+2(m+2).$$

- For $n>1$ and $m>1$, the numbers $N=2^{n+1}\cdot(2m+1)$ can be written as the sum of four successive odd positive numbers:

$$[(N/4)-3]+[(N/4)-1]+[(N/4)+1]+[(N/4)+3].$$

Summarizing,

**Proposition 1.9** *Every even number of the form $2^{n+1}\cdot(2m+1)$, where $m$ and $n$ are positive integers and $m>1$, can be written either as the sum of four successive even positive numbers if $n=1$, or as the sum of four successive odd positive numbers if $n>1$.*

In particular,

**Proposition** 1.9′ *Every even number of the form $2^2\,p$, where $p$ is prime $>3$, can be expressed as the sum of four successive even positive numbers*

$$(p-3)+(p-1)+(p+1)+(p+3).$$

*For $p=3$, we have $2^2\cdot 3=2+4+6$.*

Propositions 1.8′ and 1.9′ complement each other in that the numbers $2^2\cdot p$ have no other even or odd series representations except those stated.

We conclude the inquiry of Problem 1 by noting the key results. To represent a natural number as a sum of successive odd and/or even positive numbers is conditioned by its prime factorization form. When a number consists entirely of even or entirely of odd prime factors, then it is representable only by sums of successive odd positive numbers [Theorem 1A]. When a number consists of a single even prime factor and of at least one odd prime factor, then it is represent-able only by sums of successive even positive numbers [Theorem 1B]. When a number consists of more than one even prime factor and of at least one odd prime factor, then it is representable by both sums [Theorem 1C]. Lastly, when a number is prime, then it is representable by neither of these sums [Theorem 1.1′]. Expressed by sums of successive odd numbers starting from $1$, are all the squared numbers [Corollary 1.1] and expressed by sums of successive even numbers starting from $2$, are the products of any two consecutive natural numbers [Corollary 1.2].

## PROBLEM 2 *It is required to represent every positive integer $\succ 1$ by sums of consecutive positive integers.*

Let $a, a+1, a+2, \ldots, b$ be a sequence of $r$ consecutive positive integers. Then any other positive integer $N$ can be expressed by the series $a + (a+1) + (a+2) + \cdots + b$ if $N = (1/2) \cdot r \cdot (r + 2a - 1)$. On the other hand, $N$ can be always written as $N = d \cdot d'$, where $d$ is any of the divisors of $N$. Hence, we must have: $2 d \cdot d' = r \cdot (r + 2a - 1)$. Since $a \geq 1$, then $r + 2a - 1 > r$. Furthermore, since $r$ and $r + 2a - 1$ have different parity, then at least one of the divisors $d, d'$ must be odd. For $d = 1$, we have that $r = 2$ and $d' = 2a + 1$ must also be odd. Therefore,

**THEOREM** 2′ *No number of the form $2^{n-1}$, where $n$ is a positive integer, is expressible as a sum of consecutive positive integers.*

For $d = 1$ and $d' = 2a + 1$, one finds $a = (d' - 1)/2$ and $b = (d' + 1)/2$. Therefore,

**Proposition 2.1** *Every odd positive number can be written as the sum of two consecutive positive integers.*

In particular, we have

**Proposition 2.**1′ *Odd primes are expressible exclusively by a sum of two consecutive positive integers. For every odd prime $p$, $p = (p-1)/2 + (p+1)/2$.*

For any given pair of complementary divisors $d, d'$ of a natural number $N$ where $d$ is odd and $> 1$, there are two possible alternatives: either $d < 2d'$ or $d > 2d'$. If $d < 2d'$, then we must have $r = d$ and $r + 2a - 1 = 2d'$, which implies that $a = d' - (d-1)/2$ and $b = d' + (d-1)/2$. If $d > 2d'$, then we must have $r = 2d'$ *and* $r + 2a - 1 = d$, which implies that $a = -d' + (d+1)/2$ and $d' + (d+1)/2$. Consequently, we may state the following Theorem:

**THEOREM 2** *Every natural number $N$ with at least one odd divisor $d > 1$, is expressible as a sum of consecutive integers. Thus, for each odd divisor $d > 1$ of $N = d \cdot d'$, $N$ may be written either as*

$$\left(d' - \frac{d-1}{2}\right) + \left(d' - \frac{d-1}{2} + 1\right) + \cdots + \left(d' + \frac{d-1}{2}\right) \tag{8}$$

*if $d < 2d'$, or as*

$$\left(\frac{d+1}{2} - d'\right) + \left(\frac{d+1}{2} - d' + 1\right) + \cdots + \left(\frac{d-1}{2} + d'\right) \tag{9}$$

*if $d > 2d'$.*

For $N = p[d = p, d' = 1]$, Theorem 2 reduces to Proposition $2.1'$. The series in (8) and (9) consist of $d$ and $2d'$ positive integers, respectively.

For $d' = (d \pm 1)/2$, $N$ has the form $d \cdot (d \pm 1)/2$ and the first term in the series representation of $N$ is $1$. Thus the following Corollary may be stated:

**Corollary 2.1** *Every natural number of the form $n \cdot (n+1)/2$ can be expressed as the sum of the first $n$ positive integers:*

$$\frac{n \cdot (n+1)}{2} = 1 + 2 + 3 + \cdots + n. \qquad (10)$$

Table 8 displays the first such numbers.

| n | Number | Sum |
|---|---|---|
| 2 | 3 | 1 + 2 |
| 3 | 6 | 1 + 2 + 3 |
| 4 | 10 | 1 + 2 + 3 + 4 |
| 5 | 15 | 1 + 2 + 3 + 4 + 5 |
| 6 | 21 | 1 + 2 + 3 + 4 + 5 + 6 |
| 7 | 28 | 1 + 2 + 3 + 4 + 5 + 6 + 7 |
| 8 | 36 | 1 + 2 + 3 + 4 + 5 + 6 + 7 + 8 |
| 9 | 45 | 1 + 2 + 3 + 4 + 5 + 6 + 7 + 8 + 9 |

**Table 8**

By Theorem 2 every natural number not of the form $2^{n-1}$ is representable as a sum of consecutive positive integers. In fact, Theorem 2 provides the framework for obtaining all such representations for any natural number.

According to whether a given natural number is representable by sums of an odd and/or even number of consecutive positive integers, Theorem 2 may be decomposed into three distinct components: Theorem 2A deals with those numbers which can be expressed exclusively by sums of an odd number of consecutive integers, Theorem 2B deals with those numbers which can be expressed exclusively by sums of an even number of consecutive integers, and Theorem 3C deals with the numbers which can be expressed by both these sums.

**THEOREM 2A** *Expressed exclusively by sums of an odd number of consecutive positive integers ere all the even numbers of the form $2^n \cdot (2m+1)$ for which $2m+1 < 2^{n+1}$ with $m$ and $n$ positive integers.*

A natural number $N$ is expressed by sums of an odd number of consecutive positive integers if $d < 2d'$ for each possible pair of complementary divisors $d$ and $d'$ with $d$ odd and $> 1$. This is true only if $N$ is even and the product of all its odd prime factors is less than twice the product of all its even prime factors. Thus, let $N$ in its prime factorization form be $2^n p_1 p_2 \cdots p_L$, where $n$ is a positive integer and $p_1, p_2, \ldots, p_L$ are odd primes not necessarily distinct. Then, if $p_1 p_2 \cdots p_L < 2 \cdot 2^n$, it is obvious that the products of any combination of these primes will also be $< 2^{n+1}$ and, consequently, all series representations of $N$ are deduced by (8) for each odd divisor $d > 1$ so that $d < 2d'$.

As a numerical illustration, we list in Table 9 the series representations of the lowest order even numbers (smallest $n$) for $m = 1, 2, 3, \ldots, 10$.

| m | n ≥ | $N = 2^n\,(2\,m+1)$ | Sums | # of terms |
|---|---|---|---|---|
| 1 | 1 | $6 = 2\cdot 3$ | $1+2+3$ | 3 |
| 2 | 2 | $20 = 2^2\cdot 5$ | $2+3+4+5+6$ | 5 |
| 3 | 2 | $28 = 2^2\cdot 7$ | $1+2+3+4+5+6+7$ | 7 |
| 4 | 3 | $72 = 2^3\cdot 3^2$ | $23+24+25$<br>$4+5+6+\cdots+12$ | 3<br>$3^2$ |
| 5 | 3 | $88 = 2^3\cdot 11$ | $3+4+5+\cdots+13$ | 11 |
| 6 | 3 | $104 = 2^3\cdot 13$ | $2+3+4+\cdots+14$ | 13 |
| 7 | 3 | $120 = 2^3\cdot 3\cdot 5$ | $39+40+41$<br>$22+23+24+\cdots+26$<br>$1+2+3+\cdots+15$ | 3<br>5<br>$3\cdot 5$ |
| 8 | 4 | $272 = 2^4\cdot 17$ | $8+9+10+\cdots+24$ | 17 |
| 9 | 4 | $304 = 2^4\cdot 19$ | $7+8+9+\cdots+25$ | 19 |
| 10 | 4 | $336 = 2^4\cdot 3\cdot 7$ | $111+112+113$<br>$45+46+47+\cdots+51$<br>$6+7+8+\cdots+26$ | 3<br>7<br>$3\cdot 7$ |

All possible series representations of a number $N = 2^n\cdot(2\,m+1)$ in terms of consecutive positive integers, depend entirely on how many divisors the factor $2\,m+1$ has. In particular, when $m$ is prime there is only one such representation for $N$. Therefore, we have the following Proposition.

**Proposition 2.2** *Even numbers of the form $2^n\,p$ where $n$ is a positive integer and $p$ is any odd prime for which $p < 2^{n+1}$, are expressible exclusively by a sum of $p$ consecutive positive integers:*

$$2^n\,p = \left(2^n - \frac{p-1}{2}\right) + \left(2^n - \frac{p-1}{2} + 1\right) + \cdots + \left(2^n + \frac{p-1}{2}\right). \qquad (11)$$

For $p = 2^{n+1} - 1$, (10) reduces to the sum of the first $p$ positive integers. Hence, we have the Proposition:

**Proposition** 2.2′ *The numbers $N = 2^n\,p$ where $n$ is a positive integer and $p$ is any odd prime for which $p + 1 = 2^{n+1}$, can be expressed as the sum of the first $p$ positive integers.*

In Table 10 we list the first five such numbers and their series representation.

| n | $p=2^{n+1}-1$ | $N=2^n\,p$ | Sum |
|---|---|---|---|
| 1 | 3 | 6 | $1+2+3$ |
| 2 | 7 | 28 | $1+2+3+\cdots+7$ |
| 4 | 31 | 496 | $1+2+3+\cdots+31$ |
| 6 | 127 | 8128 | $1+2+3+\cdots+127$ |
| 12 | 8191 | 33 550 336 | $1+2+3+\cdots+8191$ |

**Table 10**

It is easy to show that the numbers $N = 2^n\,p$ of Proposition $2.2'$ are perfect, i.e., they are equal to the sum of their proper divisors. Indeed, this sum is

$$\left(1+2+2^2+\cdots+2^n\right)\cdot(p+1) - N,$$

which equals

$$\left(2^{n+1} - 1\right) \cdot 2^{n+1} - N = 2N - N = N.$$

On the other hand, the numbers $N = 2^n p$ of Proposition 2.2 for which $p < 2^{n+1} - 1$, are abundant, i.e., they are less than the sum of their proper divisors. In fact, the latter sum turns out to be $N + \left(2^{n+1} - p - 1\right)$ which always exceeds $N$. Abundant are also the numbers $2^n \cdot (2m + 1)$ where $2m + 1$ is not prime.

The series representation of any number $2^n \cdot (2m + 1)$ for which $2m + 1 < 2^n$, has the maximum number of terms for the pair of complementary divisors $d = 2m + 1$, $d' = 2^n$. Thus we may state the following Proposition:

**Proposition** 2.3 *Every number of the form $2^n \cdot (2m + 1)$ for which $2m + 1 < 2^{n+1}$ and $m, n$ are positive integers, can be written as a sum of $2m + 1$ consecutive positive integers as follows:*

$$2^n \cdot (2m + 1) = (2^n - m) + (2^n - m + 1) + \cdots + (2^n + m). \tag{12}$$

It is clear that Proposition 2.3 reduces to Proposition 2.2 if $2m + 1$ is prime.

The first term in (11) becomes $1$ for $m = 2^n - 1$. Hence we have:

**Proposition** 2.3′ *Every number of the form $2^n \cdot \left(2^{n+1} - 1\right)$ where n is any positive integer, can be written as the sum of the first $2^{n+1} - 1$ positive integers.*

If $2^{n+1} - 1$ is prime, Proposition 2.3′ reduces to Proposition 2.2′.

In Table 11 we list the first five numbers of the form $2^n \cdot \left(2^{n+1} - 1\right)$ where $2^{n+1} - 1$ is composite as well as their representations as sums of consecutive positive integers.

| n | $N = 2^n (2^{n+1} - 1)$ | First Term | Last Term | Number of Terms |
|---|---|---|---|---|
| 3 | 120 | 39<br>22<br>1 | 41<br>26<br>15 | 3<br>5<br>$3 \cdot 5$ |
| 5 | 2016 | 671<br>285<br>220<br>86<br>1 | 673<br>291<br>228<br>106<br>63 | 3<br>7<br>$3^2$<br>$3 \cdot 7$<br>$3^2 \cdot 7$ |
| 7 | 32 640 | 10 879<br>6526<br>2169<br>1912<br>615<br>342<br>1 | 10 881<br>6530<br>2183<br>1928<br>665<br>426<br>255 | 3<br>5<br>$3 \cdot 5$<br>17<br>$3 \cdot 17$<br>$5 \cdot 17$<br>$3 \cdot 5 \cdot 17$ |
| 8 | 130 816 | 18 685<br>1756<br>1 | 18 691<br>1828<br>511 | 7<br>73<br>$7 \cdot 73$ |
| 9 | 523 776 | 174 591<br>47 611<br>16 881<br>15 872<br>5586<br>1536<br>1 | 174 593<br>47 621<br>16 911<br>15 888<br>5678<br>1706<br>1023 | 3<br>11<br>31<br>$3 \cdot 11$<br>$3 \cdot 31$<br>$11 \cdot 31$<br>$3 \cdot 11 \cdot 31$ |

**Table 11**

**THEOREM 2B** *Expressed exclusively by a single sum of an even number of consecutive positive integers are the even numbers of the form $2^{n-1} p$ where $n$ is any positive integer and $p$ is any odd prime $> 2^n$. The unique series representation of these numbers is*

$$2^{n-1} p = \left(\frac{p+1}{2} - 2^{n-1}\right) + \left(\frac{p+1}{2} - 2^{n-1} + 1\right) + \cdots + \left(\frac{p-1}{2} + 2^{n-1}\right), \qquad (13)$$

*where the sum consists of $2^n$ terms. No other numbers are exclusively representable by sums of an even number of consecutive positive integers.*

A natural number $N$ is expressed exclusively by sums of an even number of consecutive positive integers if $d > 2\,d'$ for every possible pair of complementary divisors $d$ and $d'$ where $d$ is odd and $> 1$. Let $N$ in its prime factorization form be $2^m p_1 p_2 \cdots p_L$, where $m$ and $L$ are positive integers and $p_1, p_2, \cdots, p_L$ are odd primes not necessarily distinct. Then if $d > 2\,d'$ for every pair of divisors $d$ and $d'$, the integer $L$ must be $1$, i.e., N contains exactly one odd prime odd factor. For, let us suppose that integer $L$ were $2$, in which case $N$ would be of the form $2^m p_1 p_2$. Then, there would be two pairs of complementary divisors $d = p_1, d' = 2^m p_2$ and $d = p_2, d' = 2^m p_1$ satisfying simultaneously the conditions $p_1 > 2^{m+1} p_2$, $p_2 > 2^{m+1} p_1$, which is impossible. Therefore, $N$ must be of the form $2^m p$ with $p > 2^{m+1}$. Furthermore, since every odd prime is itself expressible as a sum of two consecutive positive integers (Proposition 2.1′), the integer $m$ can also be zero. Subsequently, all the numbers that are representable exclusively by sums of an even number of consecutive positive integers, must be of the form $2^{n-1}\,p$ with $p > 2^n$. Since every such number has only one odd divisor, it has a single even sum representation.

As a numerical illustration of Theorem 2B we list in Table 12 the lowest order numbers (corresponding to the smallest prime $p$ that exceeds $2^n$) for $n = 1, 2, \cdots, 10$ and their unique even sum representation.

| n | $p \ge 2^n + 1$ | $2^{n-1}\,p$ | First Term | Last Term | Number of Terms |
|---|---|---|---|---|---|
| 1 | 3 | 3 | 1 | 2 | $2^1$ |
| 2 | 5 | $2 \cdot 5$ | 1 | 4 | $2^2$ |
| 3 | 11 | $2^2 \cdot 11$ | 2 | 9 | $2^3$ |
| 4 | 17 | $2^3 \cdot 17$ | 1 | 16 | $2^4$ |
| 5 | 37 | $2^4$ | 3 | 34 | $2^5$ |
| 6 | 67 | $2^5 \cdot 67$ | 2 | 65 | $2^6$ |
| 7 | 131 | $2^6 \cdot 131$ | 2 | 129 | $2^7$ |
| 8 | 257 | $2^7 \cdot 257$ | 1 | 256 | $2^8$ |
| 9 | 521 | $2^8 \cdot 521$ | 5 | 516 | $2^9$ |
| 10 | 1031 | $2^9 \cdot 1031$ | 4 | 1027 | $2^{10}$ |

**Table 12**

For $p = 2^n + 1$, (12) reduces to the sum of the first $p - 1$ positive integers. Hence we may state the following Proposition.

**Proposition 2.4** *Every number of the form $2^{n-1}\,p$ where $n$ is a positive integer and $p$ is any odd prime for which $p - 1 = 2^n$, can be expressed uniquely as the sum of the first $2^n$ positive integers.*

From Table 12 we observe that the first primes of the form $p = 2^n + 1$ occur for $n = 1, 2, 4, 8$. One might therefore be tempted to conjecture that every number of the form $2^{2^{\nu-1}}$ is prime for $\nu$ any positive integer. However, this conjecture is first falsified for $n = 16$ since the generated number $4\,294\,967\,297$ is composite: $641 \cdot 6\,700\,417$.

It is interesting to note that the numbers which Theorem 2B refers, are deficient, i.e., the sum of their proper divisors is less than the numbers themselves. Indeed, the sum of proper divisors for these numbers is

$$\left(1 + 2 + 2^2 + \cdots + 2^{n-1}\right) \cdot (p + 1) - 2^{n-1}\,p = 2^{n-1}\,p - (p + 1 - 2^n),$$

Also deficient are the numbers of the form $2^{n-1}$ since

$$1 + 2 + 2^2 + \cdots + 2^{n-2} = 2^{n-1} - 1 < 2^{n-1}.$$

**THEOREM 2C** *Expressed by at least one sum of an odd number of consecutive positive integers and also by at least one sum of an even number of consecutive positive integers, are all the numbers of the form $2^{n-1} \cdot (2\,m + 1)$ where $m$ and $n$ are positive integers such that $2\,m + 1$ is composite and $> 2^n$.*

Let $2\,m + 1$ be any number $> 2^n$, prime or composite. By Theorem 2, every number of the form $2^{n-1} \cdot (2\,m + 1)$ can be written as the sum of $2^n$ consecutive positive integers for the complementary pair of divisors $d = 2\,m + 1$ and $d' = 2^{n-1}$. Consequently, all such numbers for which $2\,m + 1 > 2^n$, can be expressed by at least one sum of an even number of consecutive positive integers. On the other hand, let $2\,m + 1$ be any composite number. Then, $2m+1 = p \cdot (2\mu+1)$ where $p$ is an odd prime $>1$ and $\mu$ is any positive integer so that $p \le 2\mu+1$. Hence for $2m+1 > 2^n$, one has that $p < 2^n \cdot (2\mu+1)$, and for the complementary pair of divisors $d = p$ and $d' = 2^{n-1} \cdot (2\mu+1)$ it follows that the numbers $2^{n-1} \cdot (2m+1)$ can be written as the sum of p consecutive positive integers. Therefore, every number of the form $2^{n-1} \cdot (2m+1)$ for which $2m+1$ is composite and $>2^n$ can also be expressed by at least one sum of an odd number of consecutive integers.

As a numerical illustration of Theorem 2C we list in Table 13 some numbers of the form $2^{n-1}\,(2m+1)$ corresponding to the smallest $m$ which makes $2m+1$ composite and $>2^n$, as well as their representations as sums of consecutive positive integers.

| n | $m \ge 2^{n-1}$ | $2^{n-1}\,(2m+1)$ | First Term | Last Term | Number of Terms |
|---|---|---|---|---|---|
| 1 | $2^2$ | $3^2$ | 4<br>2 | 5<br>4 | 2<br>3 |
| 2 | $2^2$ | $2 \cdot 3^2$ | 5<br>3 | 7<br>6 | 3<br>$2^2$ |
| 3 | $2^2$ | $2^2 \cdot 3^2$ | 11<br>1 | 13<br>8 | 3<br>$2^3$ |
| 4 | $2*5$ | $2^3 \cdot 3 \cdot 7$ | 55<br>21<br>3 | 57<br>27<br>18 | 3<br>7<br>$2^4$ |
| 5 | $2^4$ | $2^4 \cdot 3 \cdot 11$ | 175<br>43<br>1 | 177<br>53<br>32 | 3<br>11<br>$2^5$ |
| 6 | $2^5$ | $2^5 \cdot 5 \cdot 13$ | 414<br>154<br>1 | 418<br>166<br>64 | 5<br>13<br>$2^6$ |
| 7 | $2^6$ | $2^6 \cdot 3 \cdot 43$ | 2751<br>171 | 2753<br>213<br>128 | 3<br>43<br>$2^7$ |
| 8 | $3*43$ | $2^7 \cdot 7 \cdot 37$ | 4733<br>878<br>1 | 4739<br>914<br>257 | 7<br>37<br>$2^8$ |
| 9 | $2^8$ | $2^8 \cdot 3^3 \cdot 19$ | 43 775<br>14 588<br>6903<br>4851<br>2276<br>683<br>1 | 43 777<br>14 596<br>6921<br>4877<br>2332<br>853<br>512 | 3<br>$3^2$<br>19<br>$3^3$<br>$3 \cdot 19$<br>$3^2 \cdot 19$<br>$2^9$ |

**Table 13**

It is clear that the representations of these numbers are as many as their odd divisors $>1$.

Of the numbers $2^{n-1} \cdot (2m+1)$ where $2m+1$ is composite and $>2^n$, those with two odd divisors are of the form $2^{n-1}\,p^2$ where $p$ is an odd prime for which $p^2 > 2^n$.These numbers have two representations as a sum of consecutive integers, one with an even number of terms and one with an odd number of terms. Hence, we may state the Proposition:

**Proposition 2.5** *Expressed by a single sum of an odd number of consecutive positive integers and by a single sum of an even number of consecutive positive integers, are the numbers of the form* $2^{n-1}\,p^2$ *where* $n$ *is a positive integer and* $p$ *is an odd prime for which* $p^2 > 2^2$. *These two sums are*

$$2^{n-1}\,p^2 = \left(2^{n-1}\,p - \frac{p-1}{2}\right) + \left(2^{n-1}\,p - \frac{p-1}{2} + 1\right) + \cdots + \left(2^{n-1}\,p + \frac{p-1}{2}\right), \qquad (14)$$

*and*

$$2^{n-1}\,p^2 = \left(\frac{p^2+1}{2} - 2^{n-1}\right) + \left(\frac{p^2+1}{2} - 2^{n-1} + 1\right) + \cdots + \left(\frac{p^2-1}{2} + 2^{n-1}\right), \qquad (15)$$

*consisting of* $p$ *and* $2^n$ *consecutive positive integers, respectively.*

As a numerical illustration of Proposition 2.5, we list in Table 14 (for the first ten odd primes $p$) these two series representations of the numbers $2^{n-1}\,p^2$ with the largest $n$ so that $2^n < p^2$.

| p | $2^{n-1}\cdot p^2$ | OddSum (# of terms) | EvenSum (# of terms) |
|---|---|---|---|
| 3 | $2^2\cdot 3^2$ | 11 + 12 + 13 (3) | 1 + ⋯ + 8 $(2^3)$ |
| 5 | $2^3\cdot 5^2$ | 38 + ⋯ + 42 (5) | 5 + ⋯ + 20 $(2^4)$ |
| 7 | $2^4\cdot 7^2$ | 109 + ⋯ + 115 (7) | 9 + ⋯ + 40 $(2^5)$ |
| 11 | $2^5\cdot 11^2$ | 347 + ⋯ + 357 (11) | 29 + ⋯ + 92 $(2^6)$ |
| 13 | $2^6\cdot 13^2$ | 826 + ⋯ + 838 (13) | 21 + ⋯ + 148 $(2^7)$ |
| 17 | $2^7\cdot 17^2$ | 2168 + ⋯2184 (17) | 17 + ⋯ + 272 $(2^8)$ |
| 19 | $2^7\cdot 19^2$ | 2423 + ⋯ + 2441 (19) | 53 + ⋯ + 308 $(2^8)$ |
| 23 | $2^8\cdot 23^2$ | 5877 + ⋯ + 5889 (23) | 9 + ⋯ + 520 $(2^9)$ |
| 29 | $2^8\cdot 29^2$ | 7410 + ⋯ + 7438 (29) | 165 + ⋯ + 676 $(2^9)$ |
| 31 | $2^8\cdot 31^2$ | 7921 + ⋯ + 7951 (31) | 225 + ⋯ + 736 $(2^9)$ |

**Table 14**

We conclude the inquiry of Problem 2 by noting the key results. All numbers with no odd prime factors cannot be expressed by sums of consecutive positive integers [Theorem 2′] . Expressed by sums of an odd number of consecutive positive integers are the even numbers for which the odd part is less than twice their even part [Theorem 2A]. All odd primes as well as the even numbers the odd part of which is a prime and is more than twice their even part, can be expressed by a single sum of an even number of consecutive positive integers [Theorem 2B]. The composite odd numbers as well as the even numbers for which the odd part is more than twice their even part, can be expressed by at least one sum of an odd number of consecutive positive integers and by at least one sum of an even number of consecutive positive integers [Theorem 2C].

# REFERENCES

1. L. E. Dickson, History of the Theory of Numbers, Vol.2, Chelsea, N.Y.1996.

2. T. L. Heath, A History of Greek Mathematics, Vol.1, Dover, N.Y.1981.

3. Diophantus, Arithmetica (On Polygonal Numbers, Prop.5), ed. P. Tannery, B. C. Teubner, Leipzig 1895; Eng. trans. by T. L. Heath, Cambridge 1910.

4. Nicomachus, Nicomachi Geraseni Pythagorei Introductionis Arithmeticae Libri II, II.11.1-4, R. Hoche, 1866; Eng. trans. Great Books of the Western World, vol 11, Encyclopedia Britannica, Univ. of Chicago, 1952.

5. D. A. Sardelis and T. M. Valahas, On Multidimensional Pythagorean Numbers, arxXiv : 0805.4070 v1.

# APPENDIX

The Table that follows is intended to illustrate that all possible (simple) arithmetic series representations of the first hundred natural numbers are indeed classified according to the six Theorems presented in this work. In the columns entitled Theorem 1A, Theorem 1B and Theorem 1C, we list all the representations of the first hundred positive integers as arithmetic series with a common difference of $2$, i.e., as sums of successive odd positive numbers and/or as sums of successive even positive numbers. In the columns entitled Theorem 2A, Theorem 2B and Theorem 2C, we list all the representations of the first hundred positive integers as arithmetic series with a common difference of $1$, i.e., as sums of consecutive positive integers. The number of terms for each series representation is enclosed by a parenthesis placed to the right of the corresponding series.

| **N** | $\Pi p_i$ | **Theorem 1 A** | **Theorem 1 B** | **Theorem 1 C** | **Theorem 2 A** | **Theorem 2 B** | **Theorem 2 C** |
|---|---|---|---|---|---|---|---|
| **3** | $3$ | | | | | $1+2\ (2)$ | |
| **4** | $2^2$ | $1+3\ (2)$ | | | | | |
| **5** | $5$ | | | | | $2+3\ (2)$ | |
| **6** | $2\cdot 3$ | | $2+4\ (2)$ | | $1+2+3\ (3)$ | | |
| **7** | $7$ | | | | | $3+4\ (2)$ | |
| **8** | $2^3$ | $3+5\ (2)$ | | | | | |
| **9** | $3^2$ | $1+3+5\ (3)$ | | | | | $4+5\ (2)$ |
| **10** | $2\cdot 5$ | | $4+6\ (2)$ | | | $1+\cdots+4\ (2^2)$ | |
| **11** | $11$ | | | | | $5+6\ (2)$ | |
| **12** | $2^2\cdot 3$ | | | $5+7\ (2)$<br>$2+4+6\ (3)$ | $3+4+5\ (3)$ | | |
| **13** | $13$ | | | | | $6+7\ (2)$ | |
| **14** | $2\cdot 7$ | | $6+8\ (2)$ | | | $2+\cdots+5\ (2^2)$ | |
| **15** | $3\cdot 5$ | $3+5+7\ (3)$ | | | | | $7+8\ (2)$<br>$4+5+6\ (3)$<br>$1+\cdots+5\ (5)$ |
| **16** | $2^4$ | $7+9\ (2)$<br>$1+\cdots+7\ (2^2)$ | | | | | |
| **17** | $17$ | | | | | $8+9\ (2)$ | |
| **18** | $2\cdot 3^2$ | | $8+10\ (2)$<br>$4+6+8\ (3)$ | | | | $5+6+7\ (3)$<br>$3+\cdots+6\ (2^2)$ |
| **19** | $19$ | | | | | $9+10\ (2)$ | |
| **20** | $2^2\cdot 5$ | | | $9+11\ (2)$<br>$2+\cdots+8\ (2^2)$ | $2+\cdots+6\ (5)$ | | |
| **21** | 3 · 7 | 5 + 7 + 9 (3) | | | | | 10 + 11 (2)<br>6 + 7 + 8 (3)<br>1 + ⋯ + 6 (2 · 3) |
| **22** | 2 · 11 | | 10 + 12 (2) | | | 4 + ⋯ + 7 ($2^2$) | |
| **23** | 23 | | | | | 11 + 12 (2) | |
| **24** | $2^3 \cdot 3$ | | | 11 + 13 (2)<br>6 + 8 + 10 (3)<br>3 + ⋯ + 9 ($2^2$) | 7 + 8 + 9 (3) | | |
| **25** | $5^2$ | 1 + ⋯ + 9 (5) | | | | | 12 + 13 (2)<br>3 + ⋯ + 7 (5) |
| **26** | 2 · 13 | | 12 + 14 (2) | | | 5 + ⋯ + 8 ($2^2$) | |
| **27** | $3^3$ | 7 + 9 + 11 (3) | | | | | 13 + 14 (2)<br>8 + 9 + 10 (3)<br>2 + ⋯ + 7 (2 · 3) |
| **28** | $2^2 \cdot 7$ | | | 13 + 15 (2)<br>4 + ⋯ + 10 ($2^2$) | 1 + ⋯ + 7 (7) | | |
| **29** | 29 | | | | | 14 + 15 (2) | |
| **30** | 2 · 3 · 5 | | 14 + 16 (2)<br>8 + 10 + 12 (3)<br>2 + ⋯ + 10 (5) | | | | 9 + 10 + 11 (3)<br>6 + ⋯ + 9 ($2^2$)<br>4 + ⋯ + 8 (5) |
| **31** | 31 | | | | | 15 + 16 (2) | |
| **32** | $2^5$ | 15 + 17 (2)<br>5 + ⋯ + 11 ($2^2$) | | | | | |
| **33** | 3 · 11 | 9 + 10 + 11 (3) | | | | | 16 + 17 (2)<br>10 + 11 + 12 (3)<br>3 + ⋯ + 8 (2 · 3) |
| **34** | 2 · 17 | | 16 + 18 (2) | | | 7 + ⋯ + 10 ($2^2$) | |
| **35** | 5 · 7 | 3 + ⋯ + 11 (5) | | | | | 17 + 18 (2)<br>5 + ⋯ + 9 (5)<br>2 + ⋯ + 8 (7) |
| **36** | $2^2 \cdot 3^2$ | | | 17 + 19 (2)<br>10 + 12 + 14 (3)<br>6 + ⋯ + 12 ($2^2$)<br>1 + ⋯ + 11 (2 · 3) | | | 11 + 12 + 13 (3)<br>1 + ⋯ + 8 ($2^3$) |
| **37** | 37 | | | | | 18 + 19 (2) | |
| **38** | 2 · 19 | | 18 + 20 (2) | | | 8 + ⋯ + 11 ($2^2$) | |
| **39** | 3 · 13 | 11 + 13 + 15 (3) | | | | | 19 + 20 (2)<br>12 + 13 + 14 (3)<br>4 + ⋯ + 9 (2 · 3) |
| **40** | $2^3 \cdot 5$ | | | 19 + 21 (2)<br>7 + ⋯ + 13 ($2^2$)<br>4 + ⋯ + 12 (5) | 6 + ⋯ + 10 (5) | | |
| **41** | 41 | | | | | 20 + 21 (2) | |
| **42** | 2 · 3 · 7 | | 20 + 22 (2)<br>12 + 14 + 16 (3)<br>2 + ⋯ + 12 (2 · 3) | | | | 13 + 14 + 15 (3)<br>9 + ⋯ + 12 ($2^2$)<br>3 + ⋯ + 9 (7) |
| **43** | 43 | | | | | 21 + 22 (2) | |
| **44** | $2^2 \cdot 11$ | | | 21 + 23 (2)<br>8 + ⋯ + 14 ($2^2$) | | 2 + ⋯ + 9 ($2^3$) | |
| **45** | $3^2 \cdot 5$ | 13 + 15 + 17 (3)<br>5 + ⋯ + 13 (5) | | | | | 22 + 23 (2)<br>14 + 15 + 16 (3)<br>7 + ⋯ + 11 (5)<br>5 + ⋯ + 10 (2 · 3)<br>1 + ⋯ + 9 ($3^2$) |
| **46** | 2 · 23 | | 33 + 24 (2) | | | 10 + ⋯ + 13 ($2^2$) | |
| **47** | 47 | | | | | 23 + 24 (2) | |
| **48** | $2^4 \cdot 3$ | | | $23+25$ (2)<br>$14+16+18$ (3)<br>$9+\cdots+15$ $(2^2)$<br>$3+\cdots+13$ $(2\cdot 3)$ | $15+16+17$ (3) | | |
| **49** | $7^2$ | $1+\cdots+13$ (7) | | | | | $24+25$ (2)<br>$4+\cdots+10$ (7) |
| **50** | $2\cdot 5^2$ | | $24+26$ (2)<br>$6+\cdots+14$ (5) | | | | $11+\cdots+14$ $(2^2)$<br>$8+\cdots+12$ (5) |
| **52** | $3\cdot 17$ | $15+17+19$ (3) | | | | | $25+26$ (2)<br>$16+17+18$ (3)<br>$6+\cdots+11$ $(2\cdot 3)$ |
| **52** | $2^2\cdot 13$ | | | $25+27$ (2)<br>$10+\cdots+16$ $(2^2)$ | | $3+\cdots+10$ $(2^3)$ | |
| **53** | 53 | | | | | $26+27$ (2) | |
| **54** | $2\cdot 3^3$ | | $26+28$ (2)<br>$16+18+20$ (3)<br>$4+\cdots+14$ $(2\cdot 3)$ | | | | $17+18+19$ (3)<br>$12+\cdots+15$ $(2^2)$<br>$2+\cdots+10$ $(3^2)$ |
| **55** | $5\cdot 11$ | $7+\cdots+15$ (5) | | | | | $27+28$ (2)<br>$9+\cdots+13$ (5)<br>$1+\cdots+10$ $(2\cdot 5)$ |
| **56** | $2^3\cdot 7$ | | | $27+29$ (2)<br>$11+\cdots+17$ $(2^2)$<br>$2+\cdots+14$ (7) | $5+\cdots+11$ (7) | | |
| **57** | $3\cdot 19$ | $17+19+21$ (3) | | | | | $28+29$ (2)<br>$18+19+20$ (3)<br>$7+\cdots+12$ $(2\cdot 3)$ |
| **58** | $2\cdot 29$ | | $28+30$ (2) | | | $13+\cdots+16$ $(2^2)$ | |
| **59** | 59 | | | | | $29+30$ (2) | |
| **60** | $2^2\cdot 3\cdot 5$ | | | $29+31$ (2)<br>$18+20+22$ (3)<br>$12+\cdots+18$ $(2^2)$<br>$8+\cdots+16$ (5)<br>$5+\cdots+15$ $(2\cdot 3)$ | | | $19+20+21$ (3)<br>$10+\cdots+14$ (5)<br>$4+\cdots+11$ $(2^3)$ |
| **61** | 61 | | | | | $30+31$ (2) | |
| **62** | $2\cdot 31$ | | $30+32$ (2) | | | $14+\cdots+17$ $(2^2)$ | |
| **63** | $3^2\cdot 7$ | $19+21+23$ (3)<br>$3+\cdots+15$ (7) | | | | | $31+32$ (2)<br>$20+21+22$ (3)<br>$8+\cdots+13$ $(2\cdot 3)$<br>$6+\cdots+12$ (7)<br>$3+\cdots+11$ $(3^2)$ |
| **64** | $2^6$ | $31+33$ (2)<br>$13+\cdots+19$ $(2^2)$<br>$1+\cdots+15$ $(2^3)$ | | | | | |
| **65** | $5\cdot 13$ | $9+\cdots+17$ (5) | | | | | $32+33$ (2)<br>$11+\cdots+15$ (5)<br>$2+\cdots+11$ $(2\cdot 5)$ |
| **66** | $2\cdot 3\cdot 11$ | | $32+34$ (2) | | | | $21+22+23$ (3)<br>$15+\cdots+18$ $(2^2)$<br>$1+\cdots+11$ (11) |
| **67** | 67 | | | | | $33+34$ (2) | |
| **68** | $2^2\cdot 17$ | | | $33+35$ (2)<br>$14+\cdots+20$ $(2^2)$ | | $5+\cdots+12$ $(2^3)$ | |
| **69** | $3\cdot 23$ | $21+23+25$ (3) | | | | | $34+35$ (2)<br>$32+33+34$ (3)<br>$9+\cdots+14$ $(2\cdot 3)$ |
| **70** | $2\cdot 5\cdot 7$ | | $34+36$ (2)<br>$10+\cdots+18$ (5)<br>$4+\cdots+16$ (7) | | | | $16+\cdots+19$ $(2^2)$<br>$12+\cdots+16$ (5)<br>$7+\cdots+13$ (7) |
| **71** | 71 | | | | | $35+36$ (2) | |
| **72** | $2^3 \cdot 3^2$ | | | $35+37\ (2)$<br>$22+24+26\ (3)$<br>$15+\cdots+21\ (2^2)$<br>$7+\cdots+17\ (2\cdot 3)$<br>$2+\cdots+16\ (2^3)$ | $23+24+25\ (3)$<br>$4+\cdots+12\ (3^2)$ | | |
| **73** | 73 | | | | | $36+37\ (2)$ | |
| **74** | $2\cdot 37$ | | $36+38\ (2)$ | | | $17+\cdots+20\ (2^2)$ | |
| **75** | $3\cdot 5^2$ | $23+25+27\ (3)$<br>$11+\cdots+19\ (5)$ | | | | | $37+38\ (2)$<br>$24+25+26\ (3)$<br>$13+\cdots+17\ (5)$<br>$10+\cdots+15\ (2\cdot 3)$<br>$3+\cdots+12\ (2\cdot 5)$ |
| **76** | $2^2\cdot 19$ | | | $37+39\ (2)$<br>$16+\cdots+22\ (2^2)$ | | $6+\cdots+13\ (2^3)$ | |
| **77** | $7\cdot 11$ | $5+\cdots+17\ (7)$ | | | | | $38+39\ (2)$<br>$8+\cdots+14\ (7)$<br>$2+\cdots+12\ (11)$ |
| **78** | $2\cdot 3\cdot 13$ | | $38+40\ (2)$<br>$24+26+28\ (3)$<br>$8+\cdots+18\ (2\cdot 3)$ | | | | $25+26+27\ (3)$<br>$18+\cdots+21\ (2^2)$<br>$1+\cdots+12\ (2^2\cdot 3)$ |
| **79** | | | | | | $39+40\ (2)$ | |
| **80** | $2^4\cdot 5$ | | | $39+41\ (2)$<br>$17+\cdots+23\ (2^2)$<br>$12+\cdots+20\ (5)$<br>$3+\cdots+17\ (2^3)$ | $14+\cdots+18\ (5)$ | | |
| **81** | $3^4$ | $25+27+29\ (3)$<br>$1+\cdots+17\ (3^2)$ | | | | | $40+41\ (2)$<br>$26+27+28\ (3)$<br>$11+\cdots+16\ (2\cdot 3)$<br>$5+\cdots+13\ (3^2)$ |
| **82** | $2\cdot 41$ | | $40+42\ (2)$ | | | $19+\cdots+22\ (2^2)$ | |
| **83** | 83 | | | | | $41+42\ (2)$ | |
| **84** | $2^2\cdot 3\cdot 7$ | | | $41+43\ (2)$<br>$26+28+30\ (3)$<br>$18+\cdots+24\ (2^2)$<br>$9+\cdots+19\ (2\cdot 3)$<br>$6+\cdots+18\ (7)$ | | | $27+28+29\ (3)$<br>$9+\cdots+15\ (7)$<br>$7+\cdots+14\ (2^3)$ |
| **85** | $5\cdot 17$ | $13+\cdots+21\ (5)$ | | | | | $42+43\ (2)$<br>$15+\cdots+19\ (5)$<br>$4+\cdots+13\ (2\cdot 5)$ |
| **86** | $2\cdot 43$ | | $42+44\ (2)$ | | | $20+\cdots+23\ (2^2)$ | |
| **87** | $3\cdot 29$ | $27+29+31\ (3)$ | | | | | $43+44\ (2)$<br>$28+29+30\ (3)$<br>$12+\cdots+17\ (2\cdot 3)$ |
| **88** | $2^3\cdot 11$ | | | $43+45\ (2)$<br>$19+\cdots+25\ (2^2)$<br>$4+\cdots+18\ (2^3)$ | $3+\cdots+13\ (11)$ | | |
| **89** | 89 | | | | | $44+45\ (2)$ | |
| **90** | $2\cdot 3^2\cdot 5$ | | $44+46\ (2)$<br>$28+30+32\ (3)$<br>$14+\cdots+22\ (5)$<br>$10+\cdots+20\ (2\cdot 3)$<br>$2+\cdots+18\ (3^2)$ | | | | $29+30+31\ (3)$<br>$21+\cdots+24\ (2^2)$<br>$16+\cdots+20\ (5)$<br>$6+\cdots+14\ (3^2)$<br>$2+\cdots+13\ (2^2\cdot 3)$ |
| **91** | $7\cdot 13$ | $7+\cdots+19\ (7)$ | | | | | $45+46\ (2)$<br>$10+\cdots+16\ (7)$<br>$1+\cdots+13\ (13)$ |
| **92** | $2^2\cdot 23$ | | | $45+47\ (2)$<br>$20+\cdots+26\ (2^2)$ | | $8+\cdots+15\ (2^3)$ | |

| | | | | | | | | |
|---|---|---|---|---|---|---|---|---|
| **93** | $3 \cdot 31$ | $29 + 31 + 33\ (3)$ | | | | | | $46 + 47\ (2)$<br>$30 + 31 + 32\ (3)$<br>$13 + \cdots + 18\ (2 \cdot 3)$ |
| **94** | $2 \cdot 47$ | | $46 + 48\ (2)$ | | | $22 + \cdots + 25\ (2^2)$ | | |
| **95** | $5 \cdot 19$ | $15 + \cdots + 23\ (5)$ | | | | | | $47 + 48\ (2)$<br>$17 + \cdots + 21\ (5)$<br>$5 + \cdots + 14\ (2 \cdot 5)$ |
| **96** | $2^5 \cdot 3$ | | | $47 + 49\ (2)$<br>$30 + 32 + 34\ (3)$<br>$21 + \cdots + 27\ (2^2)$<br>$11 + \cdots + 21\ (2 \cdot 3)$<br>$5 + \cdots + 19\ (2^3)$ | $31 + 32 + 33\ (3)$ | | | |
| **97** | 97 | | | | | $48 + 49\ (2)$ | | |
| **98** | $2 \cdot 7^2$ | | $48 + 50\ (2)$<br>$8 + \cdots + 20\ (7)$ | | | | | $23 + \cdots + 26\ (2^2)$<br>$11 + \cdots + 17\ (7)$ |
| **99** | $3^2 \cdot 11$ | $31 + 33 + 35\ (3)$ | | | | | | $49 + 50\ (2)$<br>$32 + 33 + 34\ (3)$<br>$14 + \cdots + 19\ (2 \cdot 3)$<br>$7 + \cdots + 15\ (3^2)$<br>$4 + \cdots + 14\ (11)$ |
| **100** | $2^2 \cdot 5^2$ | | | $49 + 51\ (2)$<br>$22 + \cdots + 28\ (2^2)$<br>$16 + \cdots + 24\ (5)$<br>$1 + \cdots + 19\ (2 \cdot 5)$ | | | | $18 + \cdots + 22\ (5)$<br>$9 + \cdots + 16\ (2^3)$ |